# Effect of Deck Deterioration on Overall System Behavior, Resilience and Remaining Life of Composite Steel Girder Bridges

Amir Gheitasi[1] and Devin K. Harris[2]

[1] Graduate Research Assistant, Department of Civil and Environmental Engineering, University of Virginia, 351 McCormick Road, Charlottesville, VA, 22904-4742, USA, Phone: (906) 370-4557; Fax: (434) 982-2951; Email: agheitasi@virginia.edu (corresponding author)
[2] Assistant Professor, Department of Civil and Environmental Engineering, University of Virginia, 351 McCormick Road, Charlottesville, VA, 22904-4742, USA, Phone: (434) 924-6373 Fax: (434) 982-2951; Email: dharris@virginia.edu

**ABSTRACT**

During past few decades, several studies have been conducted to characterize the performance of in-service girder-type bridge superstructures under operating conditions. Few of these efforts have focused on evaluating the actual response of the bridge systems, especially beyond the elastic limit of their behavior, and correlating the impact of damage to the overall system behavior. In practice, most of the in-service bridge superstructures behave elastically under the routine daily traffic; however, existing damage and deteriorating conditions would significantly influence different aspects of the structural performance including reserve capacity, resilience and remaining service-life.

The main purpose of this study is to evaluate the response of composite steel girder bridges under the effect of subsurface delamination in the reinforced concrete deck. Commercial finite element computer software, ANSYS, was implemented to perform a nonlinear analysis on a representative single-span simply supported bridge superstructure. The system failure characteristics were captured in the numerical models by incorporating different sources of material non-linearities including cracking/crushing in the concrete and plasticity in steel components. Upon validation, non-linear behavior of the system with both intact and degraded configurations was used to evaluate the impact of integrated damage mechanism on the overall system performance. Reserve capacity of this bridge superstructure was also determined with respect to the nominal element-level design capacity. As vision to the future path, this framework can be implemented to evaluate the performance of other in-service bridges degraded under the effect of different damage scenarios, thus providing a mechanism to determine a measure of capacity, resilience and remaining service-life.

**INTRODUCTION**

The safety and condition of the national infrastructure have been at the forefront of public scrutiny in recent years (ASCE 2013). Bridges represent one of the critical infrastructure components within the transportation network, as they have a potential impact on everyday lives of the traveling public. What recently put the state of these structures under the radar of public opinion are the tragic bridge failures that occurred in Minnesota (AASHTO 2008) and Washington (Tanglao et al. 2013). Irrespective of source and cause of the failures, these events often result in loss of human life and significant economic hardship to the surrounding communities (MNDOT 2013). Although, the occurrence of bridge failures is somewhat rare and most often related to natural and man-made hazards, they are mostly caused by unforeseen events such as impacts, fires or flooding rather than deterioration conditions (Kumalasari 2003).

However, it is the evolution in condition states of aging bridges that plague the health of the national transportation system. Considering the various operating conditions, in-service bridges are subjected to temporal damage and deterioration mechanisms once they are put into service. According to the national bridge inventory (FHWA 2012), almost 11% of over 600,000 bridges in-service in the United States are categorized as structurally deficient. While it is not feasible to immediately repair all of the deficient bridges, this deteriorating condition does underscore the importance of quality inspection and performance assessment mechanisms to prioritize the repair efforts. This rating emphasizes the growing challenge with preservation and rehabilitation decisions for local and federal transportation agencies, which are usually behind schedule in keeping up with appropriate maintenance strategies for all deteriorated bridges.

Several studies have been recently focused on developing applicable mechanisms and methodologies to measure different sources of defects associated bridge structures. Among those, are the innovative inspection techniques and novel technologies (Lynch and Loh 2006, Shamim et al. 2008, Gangone 2008, Vaghefi et al. 2012 and 2013), which are widely being to identify the visible deteriorations and improve the confidence in locating internal degradation mechanisms. With recent applications of the condition assessment methodologies on the massive inventory of in-service bridges in the United States, a new series of challenges arise related to the large amount of data collected from routine biennial inspection, and the manpower and expertise required to interpret this data. Even more of a challenge is the use of this collected data for rational decision-making that is based on the behavior and performance of the bridge system. What is lacking is a fundamental understanding of the system-level behavioral characteristics and the potential impact of the identified deterioration conditions on the overall performance of the bridge superstructures. A correlation between the collected inspection data and the system-level performance would provide transportation officials with a mechanism to estimate the safety and a remaining service life of the structure, but also a tool to rationalize decisions regarding long-term preservation strategies.

The focus of this investigation is to integrate the effect of damage configurations into the measure of system performance of bridge superstructures. This paper explores the impact of delamination in concrete decks on the overall system behavior of a representative composite steel girder bridge. However, the implemented approach is generic, allowing for extrapolation across other bridge types having a variety of existing conditions and operational environments. The category of composite steel girder bridges has been selected in this study as they represent one of the most common structural types in service (NBI 2012) and in reality are in a serious state of disrepair. What is unique to this study is the approach that characterizes the system-level behavior of both intact and deteriorated bridge systems, and the rationale that correlates the indentified system-level characteristics to the element-level design capacity. This approach has the potential to incorporate features such as structural components interaction and system redundancy, which are not considered in the traditional design approach of bridge superstructures (AASHTO 2012) or the current load rating practices (AASHTO 2011).

**DELAMINATION IN REINFORCED CONCRETE DECKS**

Successful maintenance practices for bridges require knowledge of condition state and also an understanding of the impacts of the possible damage and deterioration mechanisms on the overall system performance and serviceability. The main types of deterioration that composite steel girder bridges experience have been well documented in recent years (BIRM 2012). Much of the degradation often manifests in the steel girders as corrosion and section loss, which are most commonly caused as a result of leaking from expansion joints near the supports. On the other hand, reinforced concrete decks suffer from a variety of deteriorating conditions associated with cracking due to the low tensile resistance of the concrete material. Structural cracks in concrete happen under the effect of dead and live load stresses, while the non-structural cracks can be attributed to the internal stresses due to thermal and shrinkage dimensional changes. Regardless of the type and cause of the cracks, they provide openings for chloride, moisture or salt penetration which would lead to corrosion in reinforcing rebars. The corrosion product can occupy up to 10 times the volume of the corroded steel that it replaces. This expansive action creates internal pressures resulting in longitudinal surface cracks, delaminations, and spalls (Figure 1). Although these damage mechanisms may not cause the structure to collapse, in most cases they do have a major impact on the serviceability and functionality, while marring the appearance of the structure.

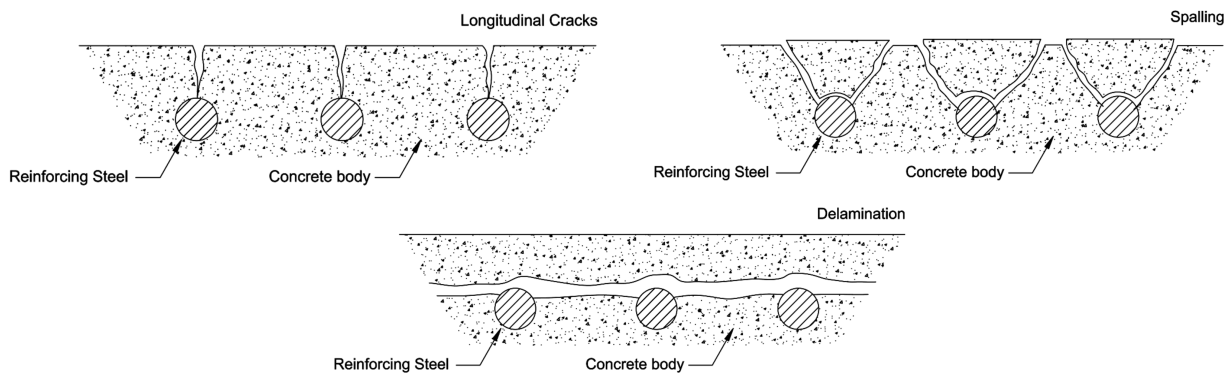

**Figure 1. Possible damage scenarios as a result of rebar corrosion.**

Delamination usually occurs as a result of separation in the concrete layers parallel and close to the surface at or near the outermost layer of reinforcing rebars (ACI 2003). Separation of the concrete layer happens when corrosion induced cracks in the vicinity of adjacent rebars join together to form a fracture plane. With a random and irregular pattern, delamination is historically considered as one of the most complicated issues associated with concrete bridges. The location and extent level of the delaminated area not only depend on the environmental conditions which directly dictates the corrosion rate, but are also related to the geometrical configuration of the concrete member such as cover thickness, rebar spacing and diameter. Once the delaminated area reaches the surface and completely separates from the concrete member, the resulting deteriorating mechanism is called spall. Delamination and spalling can occur on both the top surface and underside of an operating reinforced concrete slab (BIRM 2012).

From literature, it is evident that much of the research on delamination has focused on development of theoretical and numerical models to predict the initiation of corrosion, crack propagation, formation of the fracture plan (delamination) and local material degradation based on fracture and damage mechanics (Bažant 1979; Molina et al. 1993; Pantazopoulou and Papoulia 2001; Zhou et al. 2005; Li et al. 2006). Although the proposed models were accurate enough in predicting the behavior of deteriorated elements, they are generally too complicated to be applied to a system-level model of a bridge superstructure. This study aims to provide a simple and more practical modeling approach that can be commonly used by preservation community to evaluate the performance and serviceability of in-service bridges with delaminated concrete decks.

**METHOD OF STUDY**

To accurately evaluate the system-level behavior of in-service bridge structures under the effect of existing operating conditions, an ideal approach would be the implementation of full scale field tests on a series of representative bridges; however, this approach is neither feasible nor cost-effective. With today's computational resources and capabilities, the development of an analytical model to study the performance of intact or damaged bridge systems could be best handled numerically using a tool such as the Finite Element Method (FEM). While FEM provides an efficient mechanism to simulate the bridge system-level behavior, there are certain challenges with the modeling approach that must be properly treated to yield realistic results. Appropriate simulation assumptions, selection of material constitutive models, interpretation of complex system-level interaction and behavior, and the damage modeling techniques are some of the main challenges that should be considered in the FE modeling of bridge superstructures.

With the goal of evaluating the behavior of steel girder bridge superstructures under the effect of concrete deck delamination, the investigation approach in this study was categorized into three steps. The first step focuses on the development of system-level numerical model for a representative intact bridge superstructure. The FE model generated for this bridge system was loaded to its ultimate capacity to define the non-linear system-level behavior and characterize the failure mechanism of the system. The main outcome of this step was a fundamental understanding of the system-level performance characteristics and its correlation to the element-level behavior used in the current design practices (AASHTO 2012).

Due to the limited data that exists on the behavior of bridge structures with accumulated damage, the development of modeling strategies for integrating damage at the element-level provides a suitable alternative. As a result, the second step of this study aims at establishing a fundamental understanding on the impact of delamination on the behavior and performance of a representative reinforced concrete slab. Once the damage modeling approach within the element-level domain is validated, it can be integrated into the validated bridge model to investigate their influence on system-level behavior. A parametric investigation was also performed in this study to determine the impact of damage geometrical configurations on the system-level response, remaining service-life and susceptibility to failure. The commercial finite element (FE) computer package, ANSYS (2011), was used in this study to create f the numerical models within both element-level and system-level domains. The accuracy and validity of the FE simulation and analysis were investigated through a comparison of the numerical results to available experimental data.

**BRIDGE MODEL VALIDATION**

A full-scale laboratory investigation on a simply-supported steel girder bridge model, performed at the University of Nebraska (Kathol et al. 1995), was selected as a case study with the primary objectives of numerical simulation, analysis, and validation of the system-level behavior of composite steel girder bridges. As depicted in Figure 2a, all the structural components including the concrete deck, flexural reinforcement (in two layers), steel girder, and lateral bracing were accurately modeled with the material and geometrical properties given in the corresponding test report (Kathol et al. 1995). With no evidence of loss in the composite behavior observed during the test, the girders are assumed to be in full composite action with the concrete deck in the developed model. An ultimate capacity load test was conducted on this model bridge to evaluate the load carrying capacity of the system. Vertical concentrated loads were applied on the bridge deck to mimic two side-by-side HS-20 trucks. The corresponding loading scenario was longitudinally positioned on the system to produce the maximum positive moment at the mid-span of the bridge. According to the test setup, the model was loaded with twelve vertical concentrated loads, applied over 500 x 200 mm (20 x 8 in.) patch areas (Harris and Gheitasi 2013). Appropriate material constitutive relationships with suitable failure criteria were included in the model. Inelastic stress-strain relationships, cracking and crushing of the concrete, as well as yielding and plasticity of the steel components are the main sources of material non-linearities incorporated in the model.

Interior girder deflection at mid-span of the bridge was used to validate the proposed numerical model. As illustrated in Figure 2b, results obtained via non-linear FE analysis correlate well with experimental outcome within the elastic range of the behavior (up to 3500 kN of the applied loads). As the bridge system behave in a non-linear fashion, numerical and experimental results start to deviate with maximum error of 8%. With further application of the loads, the numerical model again begins to mimic the experimental data. Based on the validated numerical analysis, the general behavior of investigated bridge was classified into four different stages (see Figure 2b). The bridge system behaves linearly elastic prior to the formation of first flexural cracks in the concrete deck (stage A). Beyond first cracking point, it continues to carry more loads until plasticity initiates in steel girders (stage B). With further increment in the external loading, the structural stiffness drops off significantly and plastic hinges are formed in steel girders at the location of maximum moment (stage C).

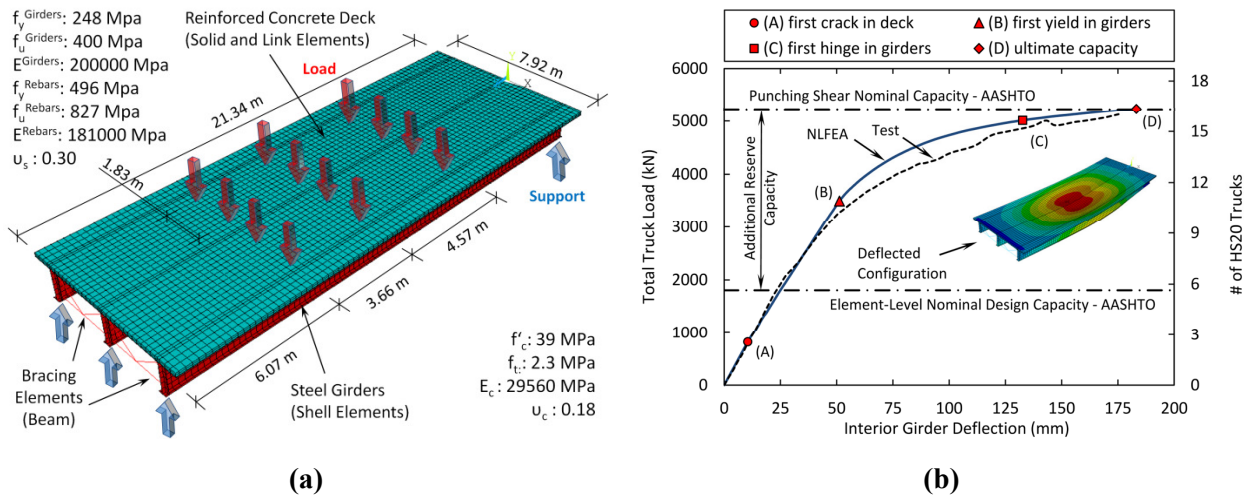

**Figure 2: Nebraska laboratory test: (a) FE model development, (b) model validation**

The laboratory experiment on this bridge model was terminated due to punching shear failure occurred under one of the patch loadings. In FE analysis, the model was loaded until it reached to the level of nominal punching shear capacity calculated based on AASHTO specifications (stage D). The ultimate load capacity and the corresponding maximum deflection predicted by the FE model were 5235 kN and 183.2 mm, respectively, as opposed to the corresponding experimental values of 5124.3 kN and 175.3 mm. In addition, Figure 3b demonstrates that the bridge system provides a significant amount of additional reserve capacity over the theoretical nominal design capacity based on AASHTO LRFD (2012) design methodology. This would validate the existence of a high level of redundancy in the system; which cannot be predicted based on the element-level design behavior of the bridge superstructure.

**DAMAGE MODELING**

Results from an experimental investigation on the behavior and performance of overlaid bridge decks, conducted at the University of California, San Diego (Seible 1998), was used in this study to develop a numerical model representing delamination. In this program, a series of simply-supported two-layer slab panels, representing an effective portion of a bridge deck spanning between longitudinal girders, were designed, cast, and tested to failure to study their interlayer behavior. A range of contact surface preparation methods (namely lubrication, roughening, and scarification with or without dowel reinforcements) between the two layers were used to provide a basis for comparison. Among all these specimens, the one that had been lubricated with bond breaking agents was selected in this study to develop the FE model. With no chemical and mechanical bonding provided between layers of concrete, this slab represents an upper bound or the worst case scenario of the magnitude of delamination that could exist in concrete decks. Figure 3a illustrates the geometrical configurations of selected slab panel. In the test setup, this slab was supported by polished and greased steel slip plates to allow horizontal translation and elastomeric bearing pads to allow rotation. A concentrated load was applied at the center of the span through a MTS actuator reacted against a structural steel load frame.

Two FE models were created in this study to assess the intact and damaged conditions of the selected slab panel. Both models were created for only one-quarter of the actual specimen to take the advantage of symmetry in geometry and loading conditions. In the numerical model of the intact slab, the concrete layers were assumed to be fully bonded at the interlayer surface, forming a monolithic configuration as desired in the test. For the damage slab model, on the other hand, the contact algorithm (ANSYS 2011) was used to simulate delamination. The contact elements only allow for compression to be transferred between the two layers; while the interlayer tension and shear forces are not restrained to mimic the situation of lubricated surfaces in the experiment. The application of contact elements would include changing-status type of non-linearity into the proposed model in addition to all other sources of material non-linearities which have been already integrated in the model.

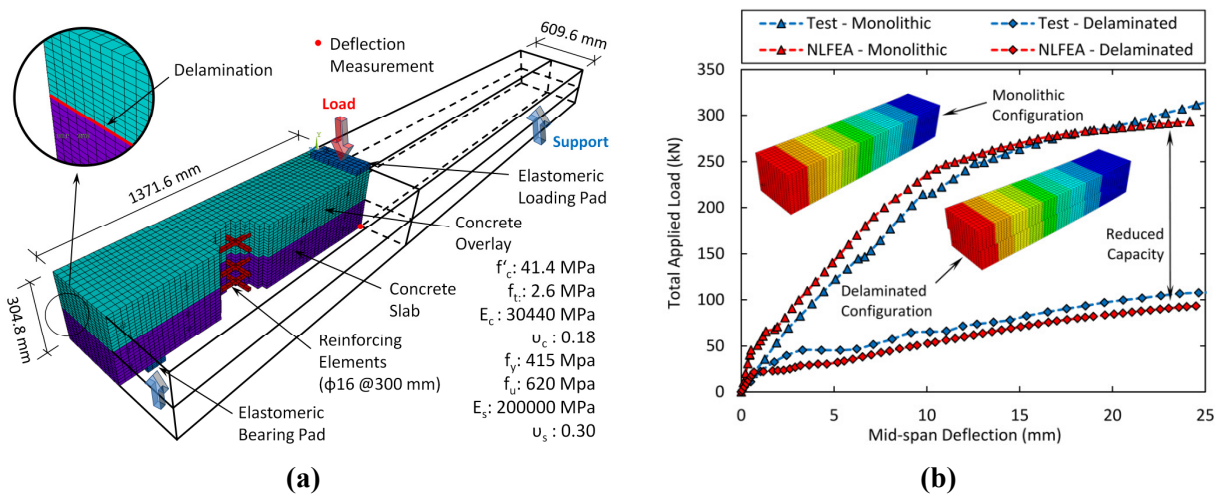

(a)          (b)

**Figure 3: Delaminated reinforced concrete slab: (a) FE model development, (b) model validation**

The load-deflection behavior at the mid-span of the simulated slabs were derived from the non-linear FE analysis and compared to the corresponding experimental result (Seible 1998) for validation purposes. Figure 3b illustrates this comparison for both the intact and deteriorated numerical models. The observed discrepancies exist between the results can be mainly attributed to the unknown material properties of the elastomeric loading pads used in the experiment. Despite the relevant numerical errors, the proposed FE models were able to predict the general non-linear behavior of both intact and damaged configurations. Moreover, comparing the results of intact and damaged models demonstrated the influence of delamination mechanism in reducing the load-carrying capacity of these slab members. This would indicate that the modeling approach implemented herein allowed for simulation of interlayer slip representing the state of delamination in concrete structures.

**DAMAGE INTEGRATION**

The validated damage modeling approach that has been used to simulate the subsurface delamination in concrete slabs can be integrated into the calibrated numerical model of the intact bridge system to identify the effect of this damage mechanism on the system-level behavior. Four different cases of delamination with respect to their level of propagation were considered in this study to update the bridge model. According to historical inspection data, delaminations in concrete decks do not usually occur in concentrated areas, but instead they are distributed randomly across the bridge deck surface area and typically occur at or above the upper layers of reinforcement. In this study, however, the delaminated area shapes and extents were assumed and positioned near midspan to induce a worst case damage configuration at this location in the simulated bridge model. As depicted in Figure 4, the assumed delaminated areas cover almost 6% to 23% of the total concrete slab area. In all cases, it was assumed that the delamination occurred at the depth of 76.2 mm (3 in.) from the top surface of the slab, as a result of corrosion in the top layer of the deck reinforcement. The fracture plane was modeled with a coupled symmetric contact elements. Only compressive stresses can be transferred towards the cracked surface, while the shear and tensile resistance at the plane of fracture have been released.

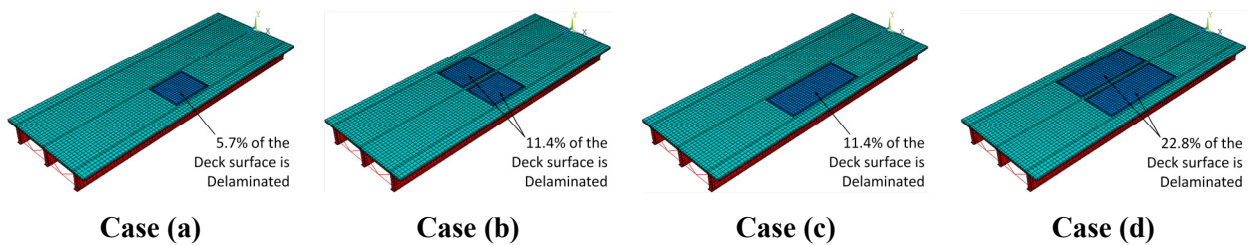

**Figure 4: Integrated damage scenarios (subsurface delamination)**

Several parameters have been proposed in the literature to measure the extent level of the concrete delamination (Li et al. 2006). Among those, the radius loss in the rebars and the crack opening in the concrete are considered as the most important parameters. Numerous experimental studies have been conducted to correlate these two parameters (Andrade et al. 1993, Alonso et al. 1998). Based on these studies, it was suggested that a 150μm reduction in the rebar radius is necessary for concrete to crack up to a width of 0.3 mm. In another study, Torres-Acosta and Martinez-Madrid (2003) found that if the crack width in the concrete slab reaches a limit within the range of 0.1-1.0 mm, delamination can occur. Based on the suggested values, the crack opening was assumed to have width of 0.75 mm in all of the analysis cases. Due to the low values of the radius reduction reported in the previous experimental investigations, the reduction in the steel rebar cross section was neglected in this study.

**DISCUSSION OF RESULTS**

Under the same loading and boundary conditions applied to the intact bridge system, the updated models were numerically analyzed to evaluate the effect of integrated damage mechanisms on the system-level behavior. Figure 5a illustrates the variations in the system-level response based on the deflection of the interior girder at mid-span of the structure. It is clear that the assumed damage configurations have negligible effects on the overall behavior and ultimate capacity of the system. In addition, the load-deflection response of the system was captured at the bottom surface of the concrete slab at the mid lateral span between the girders, see Figure 5b. Despite the local effects of the assumed damage scenarios on the deck behavior, the overall response of the deck also highlights the minor impacts of subsurface delamination on the system-level performance and maximum capacity.

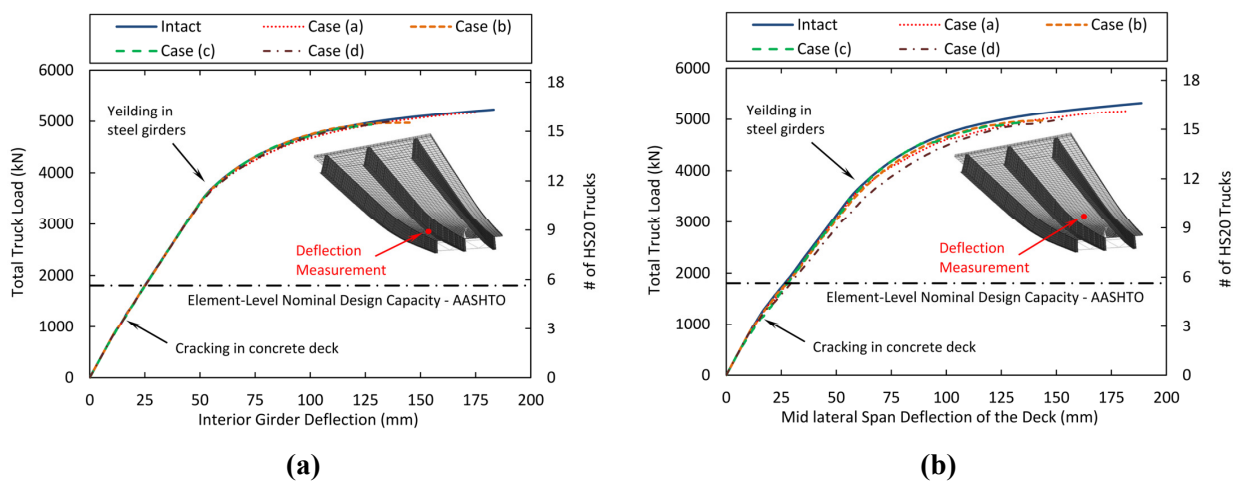

**Figure 5: Impact of damage on system-level behavior: (a) girder deflection (b) deck deflection**

As shown, all the numerical models including the intact system, demonstrated additional reserve capacity over the AASHTO LRFD element-level nominal design capacity, which would indicate the high level of system redundancy even in the presence of simulated damage mechanisms. This inherent redundancy can be attributed to the complex interaction between structural members in this type of bridge superstructures. Moreover, the behavioral stages and performance characteristics of the system, including the formation of the flexural cracks in the concrete slab and initiation of the material yielding in steel girders, appear to be unchanged for the updated models with the integrated damage scenarios. With the limited level of impact, it can be concluded that the delamination type of damage has a minimal effect on the system resiliency and remaining service life of the composite steel girder bridges. However, the assumed damage scenarios reduced the ductility of the system (ratio of maximum deflection to deflection at the first yield in girders) at the high levels of material non-linearity, due to localized failure mechanisms (concrete crushing) on the top surface of the deck in the vicinity of the delaminated regions.

Comparing Figures 5a and 5b demonstrates the fact that the overall behavior of the bridge system is highly governed by the behavior of the main load-carrying elements (i.e. girders), while the concrete deck solely distributes the lateral loads among the girders. In order to investigate the effect of subsurface delamination on the lateral load distribution mechanism in the bridge system, the principal stress vectors at the mid-span of the structure were derived from the numerical analysis for the intact and a representative damaged system (case b), and plotted in Figure 6. As demonstrated, the stress vectors in both cases have the high intensity on the top layer of the concrete slab as well as in the surrounding area of the girders, while they have a low intensity at the bottom layer of the deck between girders. This would indicate the existence of an arching action that governs the lateral load distribution mechanism in the slab on girder bridges. Moreover, the integrated damage scenario has no major effects on the load distributing mechanism of the deck system. This can be attributed to the fact that the fracture plane under the applied external loadings is subjected to compressive stresses, which are able to be fully transferred among the implemented contact elements.

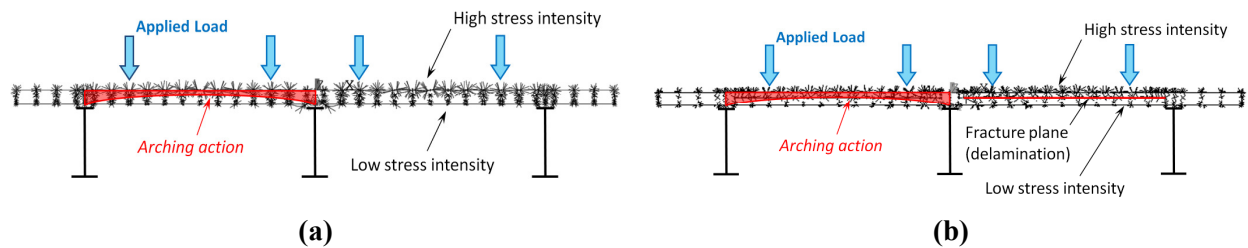

**Figure 6: Lateral load distribution mechanism: (a) intact system, (b) damaged system (case b)**

The last part of this study focuses on the impact of delamination on the deck behavior itself and the correlation of this impact to the nominal design capacities of the deck system based on the current design methodologies. According to the AASHTO LRFD (2012), the concrete decks can be designed using either Empirical or Traditional (Strip) method. In these methods, the lateral loads are assumed to be transferred among girders via arching action or flexural behavior of the deck, respectively. From the design perspective, it is assumed that the slab is vertically supported at the location of the girders, which would neglect the deflection of girders along the length of the structure. In order to address this main design assumption, the steel girders were extracted from the developed numerical models with intact and damaged configurations. Instead of girders, the slabs were supported at the location of girders through the length of structure (see Figure 7a). Under the same loading configuration, the developed models were analyzed. The load-deflection curves were derived at the mid lateral span of the deck and plotted in Figure 7b. As illustrated, the capacity of the deck decreased consistently with the expansion of the delaminated area. Regardless of the design methodology, it is clear that the actual capacity of the deck system is higher than the assumed nominal design capacities. This can be attributed to the plate action of the deck in transferring loads to the longitudinal supports and the post-cracking behavior of concrete material, which are not considered in the current design practices.

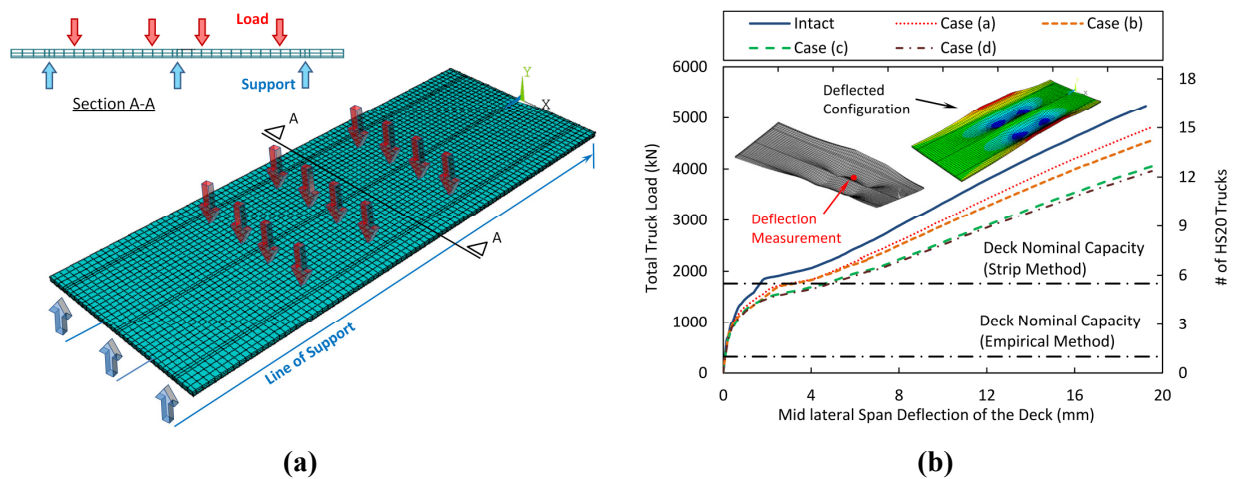

          (a)           (b)

**Figure 7: Concrete deck design capacity: (a) model development, (b) influence of damage**

## CONCLUSIONS

The main objective of this investigation was to study the influence of subsurface delamination in concrete slabs on the overall performance and behavior of composite steel girder bridge superstructures. A representative bridge model which has been previously tested to failure was selected for study. The investigation was exclusively aimed at representing a conceptual schematic of a computational modeling strategy to simulate the delamination type of damage; however, this same methodology could be extended to other types of bridge superstructures with different sources of damage mechanisms (Gheitasi and Harris 2014). The described investigation approach could be beneficial to the preservation community as a mechanism to make decisions based on in-service condition and estimate the remaining service life of the existing bridges, but also has implications in the design community where a system-level design strategy would have a major impact on design economy as compared to current element-level design strategies. Moreover, the numerical modeling approach implemented in this study also has the potential to explore the implication of advances in material, design methodologies and construction practices on the long-term performance of bridge superstructures.

## REFERENCES

ACI Committee 116R-00 (2003). Cement and Concrete Terminology, Manual of Concrete Practice, Part 1; American Concrete Institute, P.O. Box 9094, Farmington Hills, MI.

Alonso, C., Andrade, C., Rodriguez, J., and Diez, J. M. (1998). Factors controlling cracking of concrete affected by reinforcement corrosion. Mater. Struct., 31(211), 435‑441.

American Association of State and Highway Transportation Officials (2008). Bridging the gap: restoring and rebuilding the nation's bridges, Washington, D.C.

American Association of State Highway and Transportation Officials (2011). The manual for bridge evaluation. 2nd ed., American Association of State Highway and Transportation Officials, Washington, D.C.

American Association of State Highway and Transportation Officials (2012). AASHTO LRFD bridge design specifications. 6th ed., American Association of State Highway and Transportation Officials, Washington, D.C.

American Society of Civil Engineers (2013). Report card for America's infrastructure.

http://www.infrastructurereportcard.org. (Accessed December 12, 2013).

Andrade, C., Molina, F. J., and Alonso, C. (1993). Cover cracking as a function of rebar corrosion: Part I - Experiment test. Mater. Struct., 26(162), 453‑464.

ANSYS (2011). User's Manual Revision 14.0. ANSYS, Inc., Canonsburg, PA, USA.

Bažant, Z. P. (1979). Physical model for steel corrosion in concrete sea structures - Application. Journal of the Structural Division, 105(6), 1155‑1166.

Federal Highway Administration, FHWA (2012). National bridge inventory database (NBI). Federal Highway Administration, Washington, D.C.

Federal Highway Administration, FHWA (2012). Bridge Inspector's Reference Manual (BIRM). Washington, D.C.

Gangone, M. V., Whelan, M. J., Janoyan, K. D., and Jha, R. (2008). Field deployment of a dense wireless sensor network for condition assessment of a bridge superstructure. SPIE, 69330.

Gheitasi, A., and Harris, D.K. (2014). A performance-based framework for bridge preservation based on damage-integrated system-level behavior. Transportation Research Board (TRB), 93nd Annual Meeting, Transportation Research Board, Washington, D.C.

Harris, D.K., Gheitasi, A. (2013). Implementation of an energy-based stiffened plate formulation for lateral load distribution characteristics of girder-type bridges. J. Eng. Struct. 54, 168-179.

Kathol, S. Azizinamini, A. and Luedke, J. (1995). Final report: strength capacity of steel girder bridges. Nebraska Department of Roads (NDOR), RES1(0099) P469.

Kumalasari, W., and Fabian, C. H. (2003). Analysis of recent bridge failures in the United States. Journal of Performance of Constructed Facilities, 17(3), 144-150.

Li, C.Q.; Zheng, J.J.; Lawanwisut, W.; and Melchers, R.E. (2006). Concrete Delamination Caused by Steel Reinforcement Corrosion. Journal of Materials in Civil Engineering, 19(7), 591-600.

Lynch, J. P., and Loh, K. J. (2006). A summary review of wireless sensors and sensor networks for structural health monitoring. Shock and Vibration Digest, 38(2), 91-128.

Minnesota Department of Transportation. Economic impacts of the I-35W bridge collapse. http://www.dot.state.mn.us/i35wbridge/rebuild/municipal-consent/economic-impact.pdf. Accessed August18, 2013.

Molina, F. J., Alonso, C., and Andrade, C. (1993). Cover cracking as a function of rebar corrosion: Part 2−Numerical model. Materials and Structures, 26(163), 532-548.

Pantazopoulou, S. J., and Papoulia, K. D. (2001). Modeling cover-cracking due to reinforcement corrosion in RC structures. Journal of Engineering mechanics, 127(4), 342‑351.

Seible, F., Latham, C., and Krishnan, K. (1998). Structural concrete overlays in bridge deck rehabilitations – experimental program. Department of Applied Mechanics and Engineering Sciences, University of California, San Diego, La Jolla, CA.

Shamim, N. P., Gregory, L. F., Sukun, K., and David, E. C. (2008). Design and Implementation of Scalable Wireless Sensor Network for Structural Monitoring. Journal of Infrastructure Systems, 14(1), 89-101.

Tanglao, L., James, M. S., and Castellano, A. (2013). 3 Injured after I-5 bridge collapse sends cars into water, ABC News, http://abcnews.go.com. Accessed August 18, 2013.

Torres-Acosta, A. A., and Martinez-Madrid, M. M. (2003). Residual life of corroding reinforced concrete structures in marine environment. J. Mater. Civ. Eng., 15(4), 344‑353.

Vaghefi, K., Oats, R., Harris, D., Ahlborn, T., Brooks, C., Endsley, K., Roussi, C., Shuchman, R., Burns, J., and Dobson, R. (2012). Evaluation of Commercially Available Remote Sensors for Highway Bridge Condition Assessment. Journal of Bridge Engineering, 17(6), 886-895

Vaghefi, K., Ahlborn, T., Harris, D., and Brooks, C. (2013). Combined Imaging Technologies for Concrete Bridge Deck Condition Assessment. J. perform. constr. facil., In Press.

Zhou, K., Martin-Pérez, B., and Lounis, Z. (2005). Finite element analysis of corrosion-induced cracking, spalling and delamination of RC bridge decks. 1st Canadian Conference on Effective Design of Structures, Hamilton, Ont., July 10-13, 187-196.